\documentclass[10pt,a4paper]{article}
\usepackage{graphics, graphicx}     

\usepackage{cite}
\usepackage{amsmath}
\usepackage{amssymb}
\usepackage{amsthm}
\theoremstyle{plain}

\newtheorem{Theorem}{Theorem}[section] %
\newtheorem{Lemma}{Lemma}[section]

\newtheorem{Proposition}{Proposition}[section]

\newtheorem{Corollary}{Corollary}[section]

\newenvironment{Proof} 
{\par\noindent{\it Proof of}} 
{\hfill$\vspace{5mm}\scriptstyle\blacksquare$} 

\numberwithin{equation}{section} 
\numberwithin{figure}{section} 
\numberwithin{table}{section} 

\begin{document}

\setcounter{page}{1}

\markboth{M.I. Isaev}{Asymptotic behavior of the number of Eulerian orientations of graphs}

\title{Asymptotic behavior of the number of Eulerian orientations of graphs}
\date{}
\author{  Mikhail Isaev}

\maketitle
\begin{abstract}
	We consider the class of simple graphs with large algebraic connectivity (the second-smallest eigenvalue of the Laplacian matrix).
	For this class of graphs we determine the asymptotic behavior of the number of Eulerian orientations. 	
	In addition, we establish some new properties of the Laplacian matrix, as well as an estimate of a 
	conditionality of matrices with the asymptotic diagonal predominance
\end{abstract}

\section{Introduction}
$\ \ \ $
	The eulerian orientation of the graph $G$  is the orientation of its edges such that for every vertex the number of incoming edges and outgoing edges are equal. 
	We denote the number of Eulerian orientations $EO(G)$. It is easy clear that $EO(G) = 0$, if the degree of at least one vertex of $G$ is odd.	
	 Eulerian orientations of the complete graph $K_n$ are called regular tournaments. 
	
		In  \cite{Brendan1990} it is shown that for even $n\rightarrow \infty$ 
	\begin{equation}\label{Eq_1_1}
		EO(K_n) = \left(\frac{2^{n+1}}{\pi n}\right)^{(n-1)/2}
		n^{1/2}e^{-1/2} \left(1+ O(n^{-1/2+\varepsilon})\right)
	\end{equation} 
	for any $\varepsilon>0$.
	
	Undirected graphs without loops and multiple edges are called simple. 
	
The problem of counting the number of the Eulerian orienations of an undirected simple graph is
complete for the class $\# P$, see \cite{Mihail1996}. Thus this problem is difficult in terms of complexity theory. 
	
In addition, for the case of loopless  $2d$-regular graph $G$ with $n$ vertices the following estimates hold, see \cite{LasVergnas1983}, \cite{Schrijver1983}:
\begin{equation}
	2^{d}\left(\frac{(2d-1)!!}{d!}\right)^{n-1}\leq EO(G) \leq \left( \frac{(2d)!}{d!\cdot d!} \right)^{n/2}.
\end{equation} 	
An improvement of the upper bound for the regular graph case and some additional studies	in this direction
were fulfilled in \cite{LasVergnas1990}.
	
	For the simple graph $G$ we define $n\times n$ matrix $Q$ by
	\begin{equation}
		Q_{jk} = 
		\left\{
			\begin{array}{cl}
			-1, & \{v_j,v_k\}\in EG,\\
			\phantom{-}d_j,& j = k,\\
			\phantom{-}0,  & \text{ otherwise,}
		\end{array}\right.
\end{equation}
	where $n = |VG|$ and $d_j$ denotes the degree of $v_j\in VG$. The matrix $Q = Q(G)$ is called the Laplacian matrix 
	of the graph $G$. The eigenvalues
	$\lambda_1 \leq \lambda_2 \leq \ldots \leq \lambda_n$ of the
matrix $Q$ are always non-negative real numbers 
and the number of zero eigenvalues of $Q$
coincides with the number of connected components of $G$, in particular, $\lambda_1 = 0$. The eigenvalue $\lambda_2$ is called the
algebraic connectivity of the graph $G$. (For more information about the spectral properties
of the Laplace matrix see, for example, \cite{Fiedler1973} and \cite{Mohar1991}.) 
	
		According to the Kirchhoff's Matrix-Tree-Theorem, see \cite{Kirchoff1847}, we have that
\begin{equation}\label{Eq_1_3}
	t(G) = \frac{1}{n}\lambda_2\lambda_3\cdots\lambda_n,
\end{equation}
where $t(G)$ denotes the number of spanning trees of the graph $G$.

	In the present work  we generalize approach  of \cite{Brendan1990}. We determine the asymptotic behavior of the number of Eulerian orientations of simple graphs with large algebraic connectivity, see Theorem \ref{main}. 
	In Section 2 we give conventions and notations and formulate the main result. 
	
	In Section 3 we prove some basic properties of the Laplacian matrix. Some statements seem to be of independent interest, for example, we obtain an estimate of a 
	conditionality of matrices with the asymptotic diagonal predominance, see Lemma \ref{Lemma_muQ}.
	Also, we prove the following property of simple graphs with large algebraic connectivity: when you remove the vertex the number of spanning trees decreases by no more than  $cn$ times, for some $c>0$ depending only on $\lambda_2/n$, see. Corollary 3.3.   
	     
	We prove the main result in Section 4.   We express $EO(G)$ in terms of an $n$-dimensional integral using Cauchy's formula. 
	The derivation of asymptotic estimation of this integral uses three lemmas (Lemma 4.1, 4.2, 4.3), whose proofs 	are given in Section 5.
	       
\section{Main result}

Let $p \geq 1$ be a real number and $\vec{x}\in \mathbb{R}^n$. We use notation
\begin{equation}
	\left\|\vec{x}\right\|_p = \left(\sum\limits_{j=1}^{n} |x_j|^p \right)^{1/p}.
\end{equation}
For $p = \infty$ we have the maximum norm 
\begin{equation}
	\left\|\vec{x}\right\|_\infty = \max_j|x_j|.
\end{equation}
The matrix norm corresponding to the $p$-norm for vectors is
\begin{equation}
	\left\|A\right\|_p = \sup_{\vec{x}\neq 0} \frac{\left\|A\vec{x}\right\|_p}{\left\|\vec{x}\right\|_p}.
\end{equation}
If  $A$ is  the matrix of self-adjoint operator (symmetric matrix) then the norm $\|A\|_2$ is equal to
the largest module of eigenvalue of $A$ and the following inequality holds: 
\begin{equation}\label{Spec_norm}
	\left\|A\right\|_p \geq \left\|A\right\|_2.
\end{equation}
For invertible matrices one can define the condition number.
\begin{equation}
	\mu_p(A) 
	{=} \left\|A\right\|_p\cdot\left\|A^{-1}\right\|_p \geq \left\|A A^{-1}\right\|_p = 1.
\end{equation}

If $f$  is bounded both above and below by  $g$ asymptotically, we use the notation
\begin{equation}\label{Big_Theta}
	f(n) = \Theta_{k_1,k_2}\left(g(n)\right),
\end{equation}
which implies as $n\rightarrow \infty$, eventually
\begin{equation}\label{condition_Big_Theta}
	k_1|g(n)| \leq |f(n)| \leq k_2|g(n)|.
\end{equation} 
When functions $f$ and $g$ depend not only on $n$, but also on other parameters $\vec{\xi}$, we use
notation (\ref{Big_Theta}) 
meaning that condition (\ref{condition_Big_Theta}) holds uniformly for all possible values of $\vec{\xi}$.

	The main result of the present work is the following theorem.
\begin{Theorem}\label{main}
	Let $G$ be simple graph  with  $n$ vertices having even degree and the algebraic connectivity
	$\lambda_2 \geq \gamma n$ for some $\gamma > 0$. 
	Then as $n\rightarrow \infty$
	\begin{equation}\label{main_eq}
		EO(G) =\Theta_{k_1,k_2} 
		\left(
		  2^{|EG|+\frac{n}{2}} \pi^{-\frac{n}{2}} \Big/ \sqrt{t(G)}  
		\right),
	\end{equation}
	where  $t(G)$ enotes the number of spanning trees
of the graph $G$ and 
	constants $k_1,k_2 > 0$ depend only on $\sigma$.      
\end{Theorem}


\noindent
{\bf Remark 2.1.} Taking into account (\ref{Eq_1_3}), the value $t(G)$ can be represented as the principal minor 
of the Laplacian matrix $Q$.

\noindent
{\bf Remark 2.2.} For the complete graph $\lambda_1 = n$, $EK_n = \frac{n(n-1)}{2}$ and $t(K_n) = n^{n-2}$. The result of Theorem \ref{main} for this case is in agreement with asymptotic formula (\ref{Eq_1_1}).

\noindent
{\bf Remark 2.3.}  There is the result on the asymptotic behavior of Eulerian circuits  analogous to Theorem  \ref{main}, see \cite{Isaev2010}.  


\section{Some basic properties of the Laplacian matrix}
In what follows we suppose that
\begin{equation}\label{G_simple}
	G \text{ is a simple graph.}
\end{equation}
The Laplacian matrix $Q$  of the graph $G$ has the
eigenvector $[1,1,\ldots,1]^T$, corresponding to the eigenvalue $\lambda_0 = 0$.
We use notation $\hat{Q} = Q + J$, where  $J$ denotes the matrix with every entry $1$. Note that $Q$ and $\hat{Q}$ 
have the same set of
eigenvectors and eigenvalues, except for the eigenvalue corresponding to the eigenvector $[1,1,\ldots,1]^T$, which equals 
$0$ for $Q$ and $n$ for $\hat{Q}$.

Using (\ref{Spec_norm}), we get that
\begin{equation}\label{Eq_3_2}
	\lambda_{n} = ||Q||_2 \leq ||\hat{Q}||_2 \leq ||\hat{Q}||_1 = \max_{j}{\sum\limits_{k=1}^{n}} |\hat{Q}_{jk}| = n.
\end{equation}

We denote by $G_r$ the graph which arises from $G$ by removing vertices $v_1, v_2, \ldots, v_r$ and
all adjacent edges.

\begin{Lemma}\label{Lemma_connectivity}
	Let condition (\ref{G_simple}) holds for graph $G$ with $n$ vertices. Then
	\begin{equation}\label{Eq_3_3}
		\lambda_2(G) \leq \frac{n}{n-1} \min_j d_j,
	\end{equation}
	\begin{equation}\label{Eq_3_4}
		\lambda_2(G_r) \geq \lambda_2(G) - r,
	\end{equation}
	where $\lambda_1(G)$ is the algebraic connectivity of $G$ and $d_j$ is the degree of the vertex $v_j \in VG$.
\end{Lemma}
The proof of Lemma \ref{Lemma_connectivity} can be found in \cite{Fiedler1973}.

\begin{Lemma}\label{Lemma_muQ}
Let $a>0$ and $I$ be identity $n\times n$ matrix. Then for any $n$ and $n\times n$ symmetric matrix $X$ such that
	 the matrix $I+X$ is nonsingular and $|X_{ij}|\leq a/n$,	
	 \begin{equation}\label{Eq_3_5}
			\mu_2(I+X) \leq \mu_\infty(I+X) \leq  C \mu_2(I+X),
	\end{equation}
	where $C$ depends only on $a$. (does not depend on $n$)
\end{Lemma}
\begin{Proof} {\it Lemma \ref{Lemma_muQ}.}
The left-hand side of (\ref{Eq_3_5}) follows from (\ref{Spec_norm}). We order the eigenvalues of $I+X$ modulo
	\begin{equation}
		|\chi_1|\leq |\chi_2| \leq \ldots \leq |\chi_n|.
	\end{equation}
		Using (\ref{Spec_norm}), we get that
	\begin{equation}\label{Eq_3_6}
		|\chi_n| = \left\|I+X\right\|_2 \leq \left\|I+X\right\|_\infty 
		\leq \left\|I\right\|_\infty + \left\|X\right\|_\infty \leq 1 + a.
	\end{equation}

We consider $\vec{x} = (x_1,\ldots, x_n)\in \mathbb{R}^n$ such that $\left\|\vec{x}\right\|_\infty = 1$. For simplicity, we assume
that  
$x_1 = \left\|\vec{x}\right\|_\infty = 1$. We denote by $\displaystyle {\cal J} = \left\{ j \ | \ x_j > \frac{1}{2a} \right\}$.

{\bf Case 1. $\displaystyle |{\cal J}|< \frac{n}{4a}.$ } Estimating the first coordinate of $(I+X)\, \vec{x}$, we get that
\begin{equation}\label{Eq_3_7}
	\begin{aligned}
	\left\| (I+X) \vec{x}\right\|_\infty \geq x_1 - \frac{a}{n} \left( \sum\limits_{j\in {\cal J}} |x_j|  +  
	\sum\limits_{j\notin {\cal J}} |x_j| \right)
	\geq\\\geq 1 - \frac{a}{n} \left( \frac{n}{4a}\cdot 1  + n \cdot \frac{1}{2a}\right)= \frac{1}{4}\left\|\vec{x}\right\|_\infty. 
	\end{aligned}
\end{equation}

{\bf Case 2. $\displaystyle |{\cal J}|\geq \frac{n}{4a}.$ } Note that
\begin{equation}
	\sqrt{n \left\| (I+X) \vec{x}\right\|_\infty^2} \geq 
	\left\| (I+X) \vec{x}\right\|_2 \geq
	|\chi_1|\cdot\left\| \vec{x}\right\|_2 \geq 
	|\chi_1|\cdot \sqrt{|{\cal J}|\cdot \frac{1}{4a^2}}\left\|\vec{x}\right\|_\infty.
\end{equation}
Then
\begin{equation}\label{Eq_3_9}
	\left\| (I+X) \vec{x}\right\|_\infty \geq \frac{|\chi_1|}{4a^{3/2}}\left\|\vec{x}\right\|_\infty.
\end{equation}

Combining(\ref{Eq_3_7}) è (\ref{Eq_3_9}), we get that 
at least 
one of the following inequalities holds.
\begin{equation}
	\left\| (I+X)^{-1} \right\|_\infty \leq 4 \ \text{ èëè } \ \left\| (I+X)^{-1} \right\|_\infty \leq \frac{4a^{3/2}}{|\chi_1|}. 
\end{equation}
From (\ref{Eq_3_6}) we have that
\begin{equation}
	\left\|I+X\right\|_\infty \leq 1 +a.
\end{equation}
Taking into account $|\chi_n| \leq 1+a$ and $\displaystyle \mu_2(I+X) = \frac{|\chi_n|}{|\chi_1|}\geq 1$, we obtain (\ref{Eq_3_5})
\end{Proof}

The proofs of Lemma 4.1, Lemma 4.2 and Lemma 4.3 are based on the following property of the Laplacian matrix. 

\begin{Corollary}
 	Let $G$ be a simple graph with  $n$ vertices and algebraic connectivity of the graph 	
 	$\lambda_2 \geq \gamma n$ for some $\gamma > 0$. Then there is some constant 
$c_\infty>0$, depending only on $\gamma$, such that
\begin{equation}\label{normQ}
		||\hat{Q}^{-1}||_1 = ||\hat{Q}^{-1}||_\infty \leq \frac{c_{\infty}}{n}. 
\end{equation}
\end{Corollary}
\begin{Proof} {\it Corollary 3.1.}
	Using (\ref{Eq_3_3}), we get that
	\begin{equation}\label{Eq_3_13}
		 d_j \geq   \lambda_2 \frac{n-1}{n} \geq \gamma(n-1) \geq \gamma n /2.
	\end{equation}
	Taking into account (\ref{Eq_3_2}),  all eigenvalues of  $\hat{Q}$ are in the interval $[\gamma n; n]$.
	Inequality (\ref{normQ}) follows easily from the assertion of Lemma \ref{Lemma_muQ} for 
	the matrix $\Omega^T\hat{Q}\Omega$, where
	\begin{equation}
		\Omega_{jk} = 
		\left\{
			\begin{array}{cl}
			\frac{1}{\sqrt{d_j+1}},& \text{ if } j=k,\\
			0 ,& \text{ otherwise.} 
			\end{array}
		\right.
	\end{equation}
\end{Proof}

The following lemma will be applied to estimate the determinant of a matrix close to
the identity matrix $I$.
\begin{Lemma}\label{Lemma_matrix1}
	Let $X$ be an $n\times n$ matrix such that $\left\|X\right\|_2 < 1$. Then for fixed $m \geq 2$ 
	\begin{equation}
		\det(I+X) = \exp \left(  \sum\limits_{r=1}^{m-1} \frac{(-1)^{r+1}}{r}\, tr (X^r) + E_m(X)   \right),
	\end{equation}
 	where \text{tr} is the trace function and
	\begin{equation}
		|E_m(X)|\leq \frac{n}{m}\,\frac{\left\|X\right\|_2^m}{1-\left\|X\right\|_2}.
	\end{equation}
\end{Lemma}
The proof of Lemma \ref{Lemma_matrix1} is based on evaluating the trace of the matrix $\ln(I+X)$, using the representation as a convergent series. 
Lemma \ref{Lemma_matrix1} was also formulated and proved in\cite{Brendan1995}.

\begin{Lemma}\label{Lemma_Gr}
		Let $G$ be a simple graph with  $n$ vertices and algebraic connectivity of the graph 	
 	$\lambda_2 \geq \gamma n$ for some $\gamma > 0$.
	et $G_1$ be the graph which arises from $G$ by removing vertex $v_1$ and all adjacent edges. 
	 Then there is a constant $c>0$ depending only on $\gamma$ such that
	 \begin{equation}\label{Eq_3_17}
	 	\det{\hat{Q}_1} \geq \frac{\det\hat{Q}}{c n }.
	 \end{equation}
\end{Lemma}
\begin{Proof} {\it Lemma \ref{Lemma_Gr}.}  
	Note that the matrix $M_{11}$ that
results from deleting the first row and the first column of 
	$\hat{Q}$ coincides with the matrix $\hat{Q}_1$ with the exception of the diagonal elements. 
	Let $\Omega$ be a diagonal matrix such that
	\begin{equation}
		\Omega_{jj} = 
		\left\{
			\begin{array}{cl}
			1,& \text{ if } \{v_1, v_j\} \in EG,\\
			0 ,& \text{ otherwise.} 
			\end{array}
		\right.
	\end{equation}
	Define $n\times n$ matrix $X$ by 
	\begin{equation}
		X_{jk} = 
		\left\{
			\begin{array}{cl}
			\frac{1}{d_1+1},& \text{ åñëè } \{v_1, v_j\}  \notin EG, \{v_1, v_k\}  \notin EG,  \text{ and $j,k\neq 1$}\\
			0 ,& \text{ otherwise.} 
			\end{array}
		\right.
	\end{equation}
	After performing one step of the Gaussian elimination for  $\hat{Q} + \Omega+ X$, we obtain that
	\begin{equation}\label{Lemma_Gr_G1}
		\det(\hat{Q} + \Omega+ X ) = (d_1+1) \det \hat{Q}_1,
	\end{equation} 
	Using (\ref{Spec_norm}), (\ref{Eq_3_13}), we have that
	\begin{equation}
		||\Omega+ X||_2 \leq ||\Omega||_2 + ||X||_2 \leq ||\Omega||_2 + ||X||_1 \leq 1+ \frac{n}{d_1+1} \leq \frac{3}{\gamma}.
	\end{equation}
	Since the algebraic connectivity $\lambda_2\geq \gamma n$
	\begin{equation}\label{Lemma_Gr_XQ2}
		||(\Omega+X)\hat{Q}^{-1}||_2 \leq ||\Omega+X||_2 ||\hat{Q}^{-1}||_2 
		\leq \frac{3}{\gamma \lambda_2} \leq \frac{3}{ \gamma^2 n}
	\end{equation}
		Combining Lemma \ref{Lemma_matrix1} and (\ref{Lemma_Gr_XQ2}), we get that as $n\rightarrow \infty$
	\begin{equation}\label{Lemma_Gr_I+X}
		\begin{array}{r}
		\displaystyle
		\det\left({I+(\Omega+X)\hat{Q}^{-1}}\right) 
		= 
		\exp\left( tr \left((\Omega+X)\hat{Q}^{-1}\right) + E_2\left((\Omega+X)\hat{Q}^{-1}\right)\right) 
		\geq
		\\\displaystyle
		\geq 
		\exp\left( -n\frac{3}{ \gamma^2 n} + O(n^{-1})\right).
		\end{array}
	\end{equation}
		Using (\ref{Lemma_Gr_G1}) and (\ref{Lemma_Gr_I+X}) we get that as $n \rightarrow \infty$
		\begin{equation}\label{Lemma_Gr_last}
			\begin{aligned}
			(d_1+1) \det \hat{Q}_{1} =
			 \det\left({I+(\Omega+X)\hat{Q}^{-1}}\right)\det\hat{Q} 
			 \geq \\ \geq
			 \det\hat{Q} \exp\left( -3/\gamma^2 + O(n^{-1})\right).
			\end{aligned}
	\end{equation}
		Since $d_1+1 \leq n$, we obtain (\ref{Eq_3_17}).
\end{Proof}
\begin{Corollary}
Let the assumptions of Lemma \ref{Lemma_Gr} hold. 
	Let $G_r$ be the graph which arises from $G$ by removing vertices $v_1, v_2, \ldots, v_r$ and all adjacent edges. 
	 Then there is a constant $c_1>0$ depending only on $\gamma$ such that
	 	 \begin{equation}\label{Eq_3_25}
	 	\det{\hat{Q}(G_r)} \geq \frac{\det\hat{Q}(G)}{(c_1 n)^r}
	 \end{equation}
	 for any $r\leq \gamma n/2$.
\end{Corollary}
\begin{Proof} {\it Corollary 3.2.}
From (\ref{Eq_3_4}) we have that 
	\begin{equation}
		\lambda_2(G_r) \geq \gamma n - r \geq \gamma n/2.
	\end{equation}
Using $r$ times the assertion of the Lemma \ref{Lemma_Gr}, we obtain (\ref{Eq_3_25}).
\end{Proof}

According to (\ref{Eq_1_3})
\begin{equation}\label{Eq_3_27}
	t(G) = \frac{1}{n}\lambda_2\lambda_2\cdots\lambda_{n-1} = \frac{\det \hat{Q}}{n^2},
\end{equation}
then the following proposition holds. 

\begin{Corollary}	
Let the assumptions of Lemma \ref{Lemma_Gr} hold. 
	 Then there is a constant $c>0$ depending only on $\gamma$ such that	 
	 \begin{equation}
	 	t(G_1) \geq \frac{t(G)}{c n},
	 \end{equation}
	 where $t(G)$ denotes the number of spanning trees of the graph $G$.
\end{Corollary}

\begin{Lemma}\label{Lemma_struct}
	Let $a>0$ and the assumptions of Lemma \ref{Lemma_Gr} hold. Then for any set  $A \subset VG$ such, that
	  $|A| \geq an$, there is a function $h: VG \rightarrow \mathbb{N}_0$, 
	  having following properties:
	   \begin{equation}
	  	h(v) = 0, \text { åñëè } v \in A, \ \ \	h(v) \leq H,\text { äëÿ ëþáûõ }  v \in VG,
	  \end{equation}
	 	\begin{equation}\label{Eq_3_30}
	  	 \Big|\left\{w \in VG \ | \ (w,v)\in EG \text{ è } h(w)<h(v)\right\}\Big| \geq \alpha n, \text { åñëè } v \notin A, 
	  \end{equation}
	  where constants $H, \alpha > 0$ depend only on  $a$ and  $\gamma$.
\end{Lemma}
\begin{Proof} {\sc Ëåììû \ref{Lemma_struct}.}
	At first, we construct the set $A_1 = \left\{v \in VG \ | \ h(v) = 1 \right\}$, having
property (\ref{Eq_3_30}). 
	
	  If $|A| > n - \gamma n/4$, then let  $A_1 = \left\{v \in VG \ | \ v \notin A\right\}.$
	 Taking into account (\ref{Eq_3_13}), we get
that property (\ref{Eq_3_30}) hold for $\alpha = \gamma/4$. In this case $H = 1$. 
	 
	 For $|A| \leq n - \gamma n/4$ define $\vec{x} \in \mathbb{R}^n$ such that
	 \begin{equation}
	 	x_j = 
	 	\left\{
			\begin{array}{cl}
			1 - |A|/n,& \text{  }  v_j \in A,\\
			-|A|/n ,& \text{ }  v_j \notin A. 
			\end{array}
		\right.
	 \end{equation}
	 Since $x_1+x_2+\ldots+x_n =0$
	 \begin{equation}\label{Eq_3_32}
	 	\begin{aligned}
	 		\vec{x}^T Q \vec{x} = \vec{x}^T \hat{Q} \vec{x} \geq \lambda_2 \|\vec{x}\|_2^2 \geq \lambda_2 |A|\, \left(\frac{n - |A|}{n}\right)^2 
	 		\geq \\ \geq \gamma n \, a n \left(\gamma/4\right)^2 = \frac{a\gamma^3 n^2}{16}.
	 	\end{aligned}
	 \end{equation}
	  On the other hand,
	 \begin{equation}
	 	\vec{x}^T Q \vec{x} = \sum\limits_{\{v_j,v_k\}\in EG} (x_j - x_k)^2,
	 \end{equation} 
	 which is equal to the number of edges $(v,w)\in EG$, where $v\in A, w \notin A.$ 
	 We denote $A_1$ the
set of vertices $w \notin A$, having at least $\alpha n$ adjacent vertices in $A$, 
	 where $\alpha = \frac{1}{32}a\gamma^3$. 
	 \begin{equation}\label{Eq_3_34}
	 			\vec{x}^T Q \vec{x} \leq n |A_1| + \alpha n |VG|.
	 \end{equation}
	 Combining (\ref{Eq_3_32}) and (\ref{Eq_3_34}), we get that $|A_1| \geq \alpha n.$
	 
	We make further construction of the function $h$ inductively, using for the $k$-th step the
set $A^{(k)} = A\cup A_1\cup \ldots \cup A_k.$ 
	  The number of steps does not exceed $1/\alpha$ as $|A_k| \geq \alpha n$ for each
step, perhaps with the exception of the last one.
\end{Proof}
\section{Proof of Theorem \ref{main}}
In a similar way as in  \cite{Brendan1990} (see the proof of Theorem 3.1)
we note that the function
\begin{equation}
	\prod\limits_
		{\{v_j,v_k\} \in EG}
		({x_j}^{-1} x_k + {x_k}^{-1} x_j) 
\end{equation}
is the generating function the number of orientations of graph $G$ by the differences in the numbers of 
incoming and outgoing edges at each vertex.
The value $EO(G)$ is the constant term, which we can extract via Cauchy's Theorem using the unit circle as a contour for each
variable:
\begin{equation}
	EO(G) = \frac{1}{(2\pi i)^n} \oint\cdots\oint 
	\frac
		{\prod\limits_
		{\{v_j,v_k\} \in EG}
		({x_j}^{-1} x_k + {x_k}^{-1} x_j)}
		{x_1x_2\cdots x_n}dx_1 dx_2\ldots dx_n.
\end{equation} 
Making the substitution $x_j = e^{i\theta_j}$ for each $j$, we find that
\begin{equation}\label{Eq_4_2}
	EO(G) = {2^{|EG|}}{\pi^{-n}}S, \ \ \ \ \ S =\int\limits_{U_n(\pi/2)} 
     \prod\limits_{\{v_j,v_k\}\in EG} \cos\Delta_{jk} \
 d\vec{\theta},
\end{equation}
where $\Delta_{jk} = \theta_j - \theta_k$,
\begin{equation} 
U_n(\rho) = \{(x_1, x_2, \ldots, x_n)\ | \ |x_j| \leq \rho \text{ äëÿ âñåõ } j\},
\end{equation}  
and using the fact that for graphs with vertices of even degree the integrand is unchanged by the
substitutions $\theta_j \rightarrow \theta_j +  \pi$.

Let's start  the evaluation $S$ from the part that makes a major contribution to the integral
We fix some  sufficiently small constant $\varepsilon > 0$. Let
\begin{equation}
	V_0 = \{\vec{\theta}\in U_n(\pi/2) \ | \ |\theta_j - \bar{\theta}| \, (\, \mbox{mod} \, \pi) \leq n^{-1/2+\varepsilon}\text{, ãäå } \bar{\theta}=\frac{\theta_1+\ldots+\theta_n}{n}\}.
\end{equation}
By Taylor's theorem we have that for $\vec{\theta} \in V_0$  
	\begin{equation}\label{Eq_4_5}
		\prod\limits_{\{v_j,v_k\}\in EG} \cos \Delta_{jk} = 
		\exp\left(-\frac{1}{2}\sum\limits_{\{v_j,v_k\}\in EG}\Delta_{jk}^2 - 
			\frac{1}{12}\sum\limits_{\{v_j,v_k\}\in EG}\Delta_{jk}^4 + O(n^{-1+6\varepsilon})\right).
	\end{equation}
We denote by $S_0$ the contribution to $S$ in the integration over the region $V_0.$
\begin{Lemma}\label{Lemma_4_1}
	Let $G$ be a simple graph with $n$ vertices and the algebraic connectivity
	 	$\lambda_2 \geq \gamma n$ for some $\gamma > 0$.  	 
		Then for any $a,b>0$ as $n\rightarrow \infty$
	\begin{equation}\label{Eq_4_6}
		\int\limits_{V_0}	
			\exp\left(
			 - a\sum\limits_{\{v_j,v_k\}\in EG}\Delta_{jk}^2
			- b\sum\limits_{\{v_j,v_k\}\in EG} \Delta_{jk}^4 
			\right)
		 d\vec{\theta}
 = \Theta_{k_1,k_2} \left(n
		\int\limits_{\mathbb{R}^n}	
			e^{
			 - a\,\vec{\theta}^T \hat{Q} \vec{\theta}
			}
		 d\vec{\theta}\right),
	\end{equation}
	where constants $k_1,k_2>0$ depend only on $a$, $b$ and $\gamma$. 
\end{Lemma}

Lemma \ref{Lemma_4_1} follows from Lemma 8.3 of \cite{Isaev2010}. The proof is given in Section 5. 

Combining (\ref{Eq_3_27}), (\ref{Eq_4_2}), (\ref{Eq_4_5}), (\ref{Eq_4_6}) and
\begin{equation}
	\int\limits_{\mathbb{R}^n} e^{-a\,\vec{\theta}^T\hat{Q}\vec{\theta}} d \vec{\theta} =  \pi^{{n}/{2}} a^{-{n}/{2}} \Big/ \sqrt{\det\hat{Q}},
\end{equation}
we get that
\begin{equation}\label{S_0}
		S_0 =\Theta_{k_1,k_2} 
		\left(
		  2^{\frac{n-1}{2}} \pi^{\frac{n+1}{2}} \Big/ \sqrt{t(G)}  
		\right),
	\end{equation}
where constants $k_1,k_2>0$ depend only on $\gamma$.
	
Thus it remains to show that the other parts  are negligible	
One can show that
\begin{equation}\label{Eq_4_8}
	|\cos(x)|\leq \exp (-\frac{1}{2} x^2) \ \text{ for }  \ |x|\leq \frac{9}{16}\pi.
\end{equation}
	
	Divide the interval $[-\frac{1}{2}\pi,\frac{1}{2}\pi] \mbox{ mod } \pi$ into 32 equal intervals $H_0,\ldots,H_{31}$ such that $H_0 = [-\frac{1}{64}\pi,\frac{1}{64}\pi]$. For each $j$, define the region $W_j\subseteq U_n(\pi/2)$ as the set of points  $\vec{\theta} \in U_n(\pi/2)$, 
	having at least $\frac{1}{32}n$ coordinates in $H_j$. Clearly, the $W_j$'s cover $U_n(\pi/2)$ and also each $W_j$
can be mapped to $W_0$ by a uniform translation of the $\theta_j \text{ mod } \pi$. This mapping preserves the integrand of (\ref{Eq_4_2}),
and also maps $V_0$ to itself, so we have that
\begin{equation}
	\int\limits_{U_n(\pi/2)-V_0}\prod\limits_{\{v_j,v_k\}\in EG} \cos\Delta_{jk}\ d\vec{\theta}\leq 32Z,
\end{equation}
where
\begin{equation}
	Z=\int\limits_{W_0-V_0}\prod\limits_{\{v_j,v_k\}\in EG} |\cos\Delta_{jk}|\ d\vec{\theta}.
\end{equation}

We proceed by defining integrals $S_1$, $S_2$, $S_3$ in such a way that $Z$ áis obviously bounded
by their sum. We then show that each of them is negligible.
Let
\begin{equation}
	\begin{aligned}
	&V_1 =\{\vec{\theta}\in W_0 \text{ $|$ } |\theta_j|\geq \frac{1}{32}\pi \text{ for fewer than $n^\varepsilon$ values of $j$} \},\\
	&V_2 =\{\vec{\theta}\in V_1 \text{ $|$ } |\theta_j|\geq \frac{1}{16}\pi \text{ for at least one value of $j$} \}.
	\end{aligned}
\end{equation}  
Then our three integrals can be defined as
\begin{equation}
	\begin{aligned}
	&S_1 = \int\limits_{W_0-V_1}\prod\limits_{\{v_j,v_k\}\in EG} |\cos\Delta_{jk}|\ d\vec{\theta},\\
	&S_2 = \int\limits_{V_2}\prod\limits_{\{v_j,v_k\}\in EG} |\cos\Delta_{jk}|\ d\vec{\theta},\\
	&S_3 = \int\limits_{V_1-V_2-V_0}\prod\limits_{\{v_j,v_k\}\in EG} |\cos\Delta_{jk}|\ d\vec{\theta}.
	\end{aligned}
\end{equation}
	
	We begin with $S_1$. Let $h$ be the function from Lemma \ref{Lemma_struct} for the set 	$A = \{v_j \ | \ |\theta_j|\leq\frac{1}{64}\pi\}$.
	We denote $l_{min}$ such natural number that inequality 
	\begin{equation}
		|\theta_j| \geq \frac{1}{64}\pi (1 + l/H)
	\end{equation} 
	holds for at least $n^\varepsilon/H$ indices of the set $\{j \ | \ h(v_j) = l\}$. Existence of  $l_{min}$ follows from the
definition of the region $V_1$.  
	If $\theta_j$ and $\theta_k$ are such that
	\begin{equation} 
		|\theta_j|\geq \frac{1}{64}\pi (1 + l_{min}/H) \  \text{ and } \  |\theta_k|\leq \frac{1}{64}\pi (1 + (l_{min}-1)/H)
	\end{equation} 
	or vice versa, we have that
	$|\cos\Delta_{jk}|\leq \cos(\frac{1}{64}\pi/H).$ This includes at
least 
$ (\alpha n - n^{\varepsilon} ) {n^\varepsilon/H}$ edges $\{v_j,v_k\}\in EG$.  Using (\ref{Eq_3_2}) and (\ref{Eq_3_27}), we get that as $n\rightarrow \infty$
\begin{equation}\label{S_1}
	S_1\leq \pi^n\left(\cos\frac{\pi}{64H}\right)^{(\alpha n - n^{\varepsilon} ) {n^\varepsilon/H}} = 
	O\Big(\exp(-cn^{1+\varepsilon}) \Big) \,2^{\frac{n-1}{2}} \pi^{\frac{n+1}{2}}  \Big/ \sqrt{t(G)}
\end{equation}
for some constant $c > 0$, depending only on $\gamma$.

For $1\leq r\leq n^\varepsilon$ let $S_2(r)$ denote the contribution to $S_2$ 
of those $\theta \in V_2$ such that 
$|\theta_j|\geq\frac{1}{16}\pi$ for exactly $r$ values of $j$. 
If $|\theta_j|\leq\frac{1}{32}\pi$ and $|\theta_k|\geq\frac{1}{16}\pi$ 
or vice versa, we have that 
\begin{equation}\label{f_jkcos}
	|\cos\Delta_{jk}|\leq\cos\left(\frac{1}{32}\pi\right)
\end{equation}
This includes at least $r(\gamma n/2  - n^\varepsilon)$ edges $\{v_j,v_k\}\in EG$, because the
degree of any vertex of the graph $G$ is at least $\gamma n /2$, see (\ref{Eq_3_13}).
For pairs $(j,k)$ such that $|\theta_j|,|\theta_k|\leq\frac{1}{16}\pi$,   we use (\ref{Eq_4_8}).
We put $\vec{\theta'} = (\theta_1,\ldots,\theta_{n-r})$.
\begin{equation}\label{Eq_4_17}
	S_2(r)\leq \pi^r \left(\cos\frac{\pi}{32}\right)^{r(\gamma n/2 - n^\varepsilon)} 
	\sum\limits_{G_r} \int\limits_{ U_{n-r}(\pi/2) }
	\exp\left(-\frac{1}{2}\sum\limits_{\{v_j,v_k\}\in EG_r}\Delta_{jk}^2\right)
	d\vec{\theta'},
\end{equation}
where the first sum is over graphs, arises from $G$ by removing all possible sets of  $r$ vertices.

\begin{Lemma}\label{Lemma_expUn}
	Let the assumptions of Lemma \ref{Lemma_4_1} hold. Then
	\begin{equation}\label{eq_Lemma_expUn}
		\int\limits_{U_n(\pi/2)}\exp\left(-\frac{1}{2}\sum\limits_{\{v_j,v_k\}\in EG}\Delta_{jk}^2\right)d\vec{\theta}
		\leq \frac{ 2^{\frac{n-1}{2}} \pi^{\frac{n+1}{2}} n}{\sqrt{\det{\hat{Q}}}}.
	\end{equation}
\end{Lemma}

Lemma \ref{Lemma_expUn} was formulated and proved in \cite{Isaev2010} (see Lemma 6.1). The proof is given also in Section 5. 

Using Lemma  \ref{Lemma_connectivity} and Corollary 3.2, we get that 
\begin{equation}
 	\lambda_2(G_r) \geq \gamma n  -  n^\varepsilon, 
 \end{equation}	
 \begin{equation}\label{Eq_4_20}
 		\det{\hat{Q}(G_r)} \geq \frac{\det\hat{Q}}{\left(c_1 n\right)^r}.
 \end{equation}
According to Lemma \ref{Lemma_expUn} we have that
\begin{equation}\label{Eq_4_21}
 		\int\limits_{U_{n-r}(\pi/2)}\exp\left(-\frac{1}{2}\sum\limits_{\{v_j,v_k\}\in EG_{r}}\Delta_{jk}^2\right)d\vec{\theta'}
		\leq \frac{ 2^{\frac{n-r-1}{2}} \pi^{\frac{n-r+1}{2}} n}{\sqrt{\det{\hat{Q}(G_r)}}}.
 \end{equation}
Combining (\ref{Eq_4_17}) with (\ref{Eq_4_20}), (\ref{Eq_4_21}) and allowing $n^r$ for the choice of the set of $r$ vertices 
for $G_r$, we get that
\begin{equation}
  	 		S_2(r) \leq  2^{\frac{n-r-1}{2}} \pi^{\frac{n+r+1}{2}} n^{r+1} 
 	 		\left(\cos\frac{\pi}{32}\right)^{r(\gamma n/2 - n^\varepsilon)}  \frac{(c_1 n)^{r/2}}{\sqrt{\det\hat{Q}}}. 
 \end{equation} 
Using (\ref{Eq_3_2}) and summing over $0 \leq r\leq n^\varepsilon$, we find that
\begin{equation}\label{S_2}
 		S_2 = \sum\limits_{r=1}\limits^{n^\varepsilon}S_2(r) = O(c^{-n})\,2^{\frac{n-1}{2}} \pi^{\frac{n+1}{2}}  \Big/ \sqrt{t(G)}
 \end{equation}
for some constant $c>1$, depending only on $\gamma$.

	Note that $\Delta_{jk}\leq\frac{1}{8}\pi$ for $\vec{\theta} \in V_1 - V_2$, thus 
	\begin{equation}
		V_1 - V_2 \subset V_3 = \{\vec{\theta}\in U_n(\pi/2) \ | \ |\theta_j - \bar{\theta}| \, (\, \mbox{mod} \, \pi) \leq \pi/8 \},
	\end{equation}  
	where
	\begin{equation}
		\bar{\theta}=\frac{\theta_1+\ldots+\theta_n}{n}.
	\end{equation}
	Since the integrand is invariant under uniform
translation of all the
 $\theta_j$'s mod $\pi$ as well as  $V_0$ and  $V_3$ are mapped into itself, we can fix $\bar{\theta} = 0$ and multiply it by the ratio of its range $\pi$ to the length $n^{-1/2}$ of the vector $\frac{1}{n}[1,1,\ldots,1]^T$. 
	Thus we get that
	\begin{equation} 
		S_3 \leq \pi n^{1/2} \int\limits_{L\cap U_n({\pi}/{8}) \setminus V_0} 
		\prod\limits_{\{v_j,v_k\}\in EG}
	|\cos\Delta_{jk}|dL,
	\end{equation}
where $L$ denotes the orthogonal complement to the vector $[1,1,\ldots,1]^T$.
 In a similar way as in (\ref{Eq_4_17}) we find that
  	\begin{equation}\label{S44}
 		S_3\leq  	\pi n^{1/2}	\int\limits_{ L\cap U_n({\pi}/{8}) \setminus V_0}
	\exp\left(-\frac{1}{2}\sum\limits_{\{v_j,v_k\}\in EG}\Delta_{jk}^2\right)
	dL.
 	\end{equation} 
 \begin{Lemma}\label{Lemma_exp-V0}
		Let the assumptions of Lemma \ref{Lemma_4_1} hold. Then as $n\rightarrow \infty$
	\begin{equation}\label{eq_Lemma_exp-V0}
		\int\limits_{L \setminus U_n(n^{-1/2+\varepsilon})}\exp\left(-\frac{1}{2}\sum\limits_{\{v_j,v_k\}\in EG}\Delta_{jk}^2\right)dL
		=O\left(\exp(-cn^{2\varepsilon})\right) \frac{ 2^{\frac{n-1}{2}} \pi^{\frac{n-1}{2}} n^{1/2}}{\sqrt{\det{\hat{Q}}}}
	\end{equation}
	for some $c>0$, depending only on  $\gamma$.
\end{Lemma}

Lemma \ref{Lemma_exp-V0} was formulated and proved in \cite{Isaev2010} (see Lemma 6.2). The proof is given also in Section 5.  

Combining (\ref{Eq_3_27}), (\ref{eq_Lemma_exp-V0}) and (\ref{S44}),
we get that as  $n\rightarrow \infty$
\begin{equation}\label{S_4}
	S_3 = O\left(\exp(-cn^{2\varepsilon})\right) \frac{ 2^{\frac{n-1}{2}} \pi^{\frac{n+1}{2}} n}{\sqrt{\det{\hat{Q}}}}
		= O\left(\exp(-cn^{2\varepsilon})\right) \,2^{\frac{n-1}{2}} \pi^{\frac{n+1}{2}}  \Big/ \sqrt{t(G)}
\end{equation}
	for some $c>0$, depending only on $\gamma$.
	Combining (\ref{S_0}), (\ref{S_1}), (\ref{S_2}) and  (\ref{S_4}), we obtain (\ref{main_eq}) as well as the following lemma.
\begin{Lemma}\label{Lemma_S1234}
	Let the assumptions of Lemma  \ref{Lemma_4_1} hold. Then as $n\rightarrow \infty$
	\begin{equation}\label{eq_Lemma_S1234}
		S  = \left(1 + O\left(\exp(-cn^{2\varepsilon})\right)\right) S_0
	\end{equation}
	for some $c>0$, depending only on  $\gamma$.
\end{Lemma}
\section{Proofs of Lemmas 4.1 - 4.3}
In this section we always assume that the assumptions of Lemma \ref{Lemma_4_1}  hold.
Let
\begin{equation}
	\vec{\phi} = \vec{\phi}(\vec{\theta}) = [\phi_1(\vec{\theta}),\ldots, \phi_n(\vec{\theta})]^T
	= \hat{Q}\vec{\theta}. 
\end{equation}
Let $P(\vec{\theta)}$ be the orthogonal projection  $\vec{\theta}$ onto the space $L$, where
$L$ denotes the orthogonal complement to the vector $[1,1,\ldots,1]^T$. Note that
\begin{equation}\label{P_Q}
	Q\vec{\theta} = QP(\vec{\theta}).
\end{equation}
For any $a>0$, we have that: 
\begin{equation}\label{Integral_hatQ}
	\int\limits_{\mathbb{R}^n} e^{-a\,\vec{\theta}^T\hat{Q}\vec{\theta}} d \vec{\theta} =  \pi^{{n}/{2}} a^{-{n}/{2}} \Big/ \sqrt{\det\hat{Q}}
\end{equation}
and
\begin{equation}\label{Integral_hatQL}
	\int\limits_{L} e^{-a\,\vec{\theta}^T\hat{Q}\vec{\theta}} d L
	=
	\int\limits_{L} e^{-a\,\vec{\theta}^TQ\vec{\theta}} d L =  \pi^{\frac{n-1}{2}} a^{-\frac{n-1}{2}} n^{1/2} \Big/ \sqrt{\det\hat{Q}}.
\end{equation}

\begin{Proof} {\it Lemma \ref{Lemma_expUn}.}
	Note that
	\begin{equation}\label{sum->Q}
		\sum\limits_{\{v_j,v_k\}\in EG}\Delta_{jk}^2= \vec{\theta}^T Q \vec{\theta}.
	\end{equation}
	The diagonal of the cube $U_n(\pi/2)$ is equal to $\pi n^{1/2}$. Using (\ref{P_Q}), we find that
		\begin{equation}\label{Integral<Un}
		\int\limits_{U_n(\pi/2)}\exp\left(-\frac{1}{2}\sum\limits_{\{v_j,v_k\}\in EG}\Delta_{jk}^2\right)d\vec{\theta}
		\leq
		\pi n^{1/2} \int\limits_{L} e^{-\frac{1}{2}\,\vec{\theta}^T Q\vec{\theta}} d L.
	\end{equation}
	Combining (\ref{Integral_hatQL}) and (\ref{Integral<Un}), we obtain (\ref{eq_Lemma_expUn}).
\end{Proof}

Note that for some function $g_1(\vec{\theta}) = g_1(\theta_2,\ldots,\theta_{n})$ 
\begin{equation}\label{g_theta}
	\vec{\theta}^T\hat{Q}\vec{\theta} = \frac{\phi_1(\vec{\theta})^2}{d_1+1} + g_1(\vec{\theta}).
\end{equation}
We recall that (see (\ref{Eq_3_13})) 
\begin{equation}\label{gamman/2}
	\min_j d_j \geq \gamma n/2.
\end{equation}
Combining (\ref{g_theta}) and (\ref{gamman/2}), we get that as $n\rightarrow \infty$
\begin{equation}\label{5.9}
	\begin{aligned}
	\int\limits_{\mathbb{R}^n} e^{-a\,\vec{\theta}^T\hat{Q}\vec{\theta}} d \vec{\theta}
	&=
	\int\limits_{-\infty}\limits^{+\infty}\cdots \int\limits_{-\infty}\limits^{+\infty}
	e^{-a\,g_1(\theta_2,\ldots,\theta_{n})}
	\left( \int \limits_{-\infty}\limits^{+\infty} e^{-a\,\frac{\phi_1(\vec{\theta})^2}{d_1+1}} d \theta_1  \right)
	d \theta_2 \ldots d \theta_n\\
	&=
		\left(1 + O\left(\exp(-\tilde{c}n^{2\varepsilon})\right)\right)
			\int \limits_{|\phi_1(\vec{\theta})|\leq \frac{1}{2}c_\infty^{-1} n^{1/2+\varepsilon}} e^{-a\,\vec{\theta}^T\hat{Q}\vec{\theta}} d \vec{\theta}  
	\end{aligned}
\end{equation}
for some $\tilde{c}>0$, depending obly on $\gamma$ and $a$, where $c_\infty$ is the constant of (\ref{normQ}).
Combining
similar to (\ref{5.9}) expressions for $\phi_1, \ldots \phi_n$, we find that as $n\rightarrow \infty$
\begin{equation}\label{5.10}
	\int \limits_{||\vec{\phi}(\vec{\theta})||_\infty\leq \frac{1}{2}c_\infty^{-1} n^{1/2+\varepsilon}} e^{-a\,\vec{\theta}^T\hat{Q}\vec{\theta}} d \vec{\theta}
	= \left(1 + O\left(\exp(-cn^{2\varepsilon})\right)\right) \int\limits_{\mathbb{R}^n} e^{-a\,\vec{\theta}^T\hat{Q}\vec{\theta}} d \vec{\theta}
\end{equation}
for some $c>0$, depending only on $\gamma$ è $a$. Using (\ref{5.10}) and (\ref{normQ}), we get that as $n\rightarrow \infty$
\begin{equation}\label{1/2U_n}
		\int \limits_{U_n( \frac{1}{2}n^{-1/2+\varepsilon})} e^{-a\,\vec{\theta}^T\hat{Q}\vec{\theta}} d \vec{\theta}
	= \left(1 + O\left(\exp(-cn^{2\varepsilon})\right)\right) \int\limits_{\mathbb{R}^n} e^{-a\,\vec{\theta}^T\hat{Q}\vec{\theta}} d \vec{\theta}.
\end{equation}

\begin{Proof} {\it Lemma \ref{Lemma_exp-V0}.} 
	Note that 
	\begin{equation}
		||P(\vec{\theta})||_\infty = ||\vec{\theta} 
		- \bar{\theta} [1,1,\ldots,1]^T||_\infty \leq 2 ||\vec{\theta} ||_\infty,
	\end{equation}
	where
	\begin{equation}
 		\bar{\theta} = \frac{\theta_1 + \theta_2 + \ldots  \theta_n}{n}.
	\end{equation}
	Thus
	\begin{equation}\label{U_n_in_P}
		U_n( \frac{1}{2}n^{-1/2+\varepsilon}) \subset \left\{ \vec{\theta} \ | \ P(\vec{\theta}) \in U_n(n^{-1/2+\varepsilon})\right\}
	\end{equation}
	Using (\ref{P_Q}), (\ref{sum->Q}) and (\ref{U_n_in_P}), we find that
	\begin{equation}\label{LcapU_n->Un}
		\begin{aligned}
			\int\limits_{L \cap U_n(n^{-1/2+\varepsilon})}\exp\left(-\frac{1}{2}\sum\limits_{\{v_j,v_k\}\in EG}\Delta_{jk}^2\right)dL =
		 	\int\limits_{L \cap U_n(n^{-1/2+\varepsilon})} e^{-\frac{1}{2}\,\vec{\theta}^T Q\vec{\theta}} d L =\\
		 	= \int\limits_{P(\vec{\theta}) \in U_n(n^{-1/2+\varepsilon})} e^{-\frac{1}{2}\,\vec{\theta}^T \hat{Q}\vec{\theta}} d \vec{\theta}
		 	\Big/ \int\limits_{-\infty}\limits^{+\infty} e^{-\frac{1}{2}nx^2} dx 
		 	\geq
		 	\frac{n^{1/2}} {\sqrt{2 \pi}} \int\limits_{U_n(\frac{1}{2}n^{-1/2+\varepsilon})} e^{-\frac{1}{2}\,\vec{\theta}^T \hat{Q}\vec{\theta}} d \vec{\theta}.
		\end{aligned}
	\end{equation}
	Combining (\ref{Integral_hatQ}), (\ref{1/2U_n}) and (\ref{LcapU_n->Un}), we obtain (\ref{eq_Lemma_exp-V0}).
\end{Proof}

\begin{Proposition} As $n\rightarrow \infty$
	\begin{equation}\label{P51}
		\begin{aligned}
		\int\limits_{U_n(n^{-1/2+\varepsilon})} \exp\left( - a\, \vec{\theta}^T \hat{Q} \vec{\theta} - b \sum\limits_{\{v_j,v_k\}\in EG} 
		\Delta_{jk}^4 \right) &d\vec{\theta} = \\
		&= \Theta_{k_1,k_2} \left(\int\limits_{\mathbb{R}^n} e^{ - a\, \vec{\theta}^T \hat{Q} \vec{\theta} } d\vec{\theta} \right),
		\end{aligned}
	\end{equation}
	where constants $k_1,k_2>0$ depend only on $a$, $b$ and $\gamma$.
\end{Proposition}
In thó present paper we give only the scheme of the proof of Proposition 5.1.
The detailed proof can be found  in \cite{Isaev2010} (see Lemma 8.3). 

\noindent
{\it Scheme of the proof of Proposition 5.1.} 
Since $a,b>0$  
\begin{equation}\label{5.17}
		\begin{aligned}
		\int\limits_{U_n(n^{-1/2+\varepsilon})} \exp\left( - a\, \vec{\theta}^T \hat{Q} \vec{\theta} - b \sum\limits_{\{v_j,v_k\}\in EG} 
		\Delta_{jk}^4 \right) d\vec{\theta}  
		\leq \int\limits_{\mathbb{R}^n} e^{ - a\, \vec{\theta}^T \hat{Q} \vec{\theta} } d\vec{\theta}.
		\end{aligned}
	\end{equation}
	Let 
	\begin{equation}
		R_k(\vec{\theta}) = 8bn \sum\limits_{j=k}\limits^{n} \theta_j^4, \ \ \ 1\leq k\leq n.
	\end{equation}
	Using the representation of the integral as an iterated integral, one can
show that
	\begin{equation}\label{5.19}
		\begin{aligned}
		\int\limits_{U_n(n^{-1/2+\varepsilon})} \phi_j^4 e^{ - a\, \vec{\theta}^T \hat{Q} \vec{\theta} - R_k(\vec{\theta})} d\vec{\theta} 
			\leq c' n^2  \int\limits_{U_n(n^{-1/2+\varepsilon})} e^{ - a\, \vec{\theta}^T \hat{Q} \vec{\theta} - R_k(\vec{\theta})} d\vec{\theta} +\\+ 
			O\left(\exp(-cn^{2\varepsilon})\right)\int\limits_{\mathbb{R}^n} e^{ - a\, \vec{\theta}^T \hat{Q} \vec{\theta} } d\vec{\theta}
		\end{aligned}
	\end{equation} 
	for some constants $c,c'>0$, depending only on $a$, $b$ and $\gamma$. 
	
	Using	(\ref{5.19}), one can get that:
		\begin{equation}
		\begin{aligned}
		\int\limits_{U_n(n^{-1/2+\varepsilon})} \theta_j^4 e^{ - a\, \vec{\theta}^T \hat{Q} \vec{\theta} - R_k(\vec{\theta})} d\vec{\theta} 
			\leq \frac{c'_1}{ n^2}  \int\limits_{U_n(n^{-1/2+\varepsilon})} e^{ - a\, \vec{\theta}^T \hat{Q} \vec{\theta} - R_k(\vec{\theta})} d\vec{\theta} +\\+ 
			O\left(\exp(-c_1 n^{2\varepsilon})\right)\int\limits_{\mathbb{R}^n} e^{ - a\, \vec{\theta}^T \hat{Q} \vec{\theta} } d\vec{\theta};
		\end{aligned}
	\end{equation} 
	
	\begin{equation}\label{5.21}
		\begin{aligned}
		\int\limits_{U_n(n^{-1/2+\varepsilon})}  e^{ - a\,\vec{\theta}^T \hat{Q} \vec{\theta} - R_k(\vec{\theta})} d\vec{\theta} 
			\geq \left(1 +  \frac{c'_2}{n} \right)  \int\limits_{U_n(n^{-1/2+\varepsilon})} e^{ - a\, \vec{\theta}^T \hat{Q} \vec{\theta} - R_{k+1}(\vec{\theta})} d\vec{\theta} +\\+ 
			O\left(\exp(-c_2 n^{2\varepsilon})\right)\int\limits_{\mathbb{R}^n} e^{ - a\vec{\theta}^T \hat{Q} \vec{\theta} } d\vec{\theta},
		\end{aligned}
	\end{equation} 
	for some constants $c_1,c'_1,c_2,c'_2>0$, depending only on $a$, $b$ and $\gamma$. Note that	
	\begin{equation}\label{5.22}
		b\sum\limits_{\{v_j,v_k\}\in EG} 
		\Delta_{jk}^4  \leq 8bn \sum\limits_{j=1}\limits^{n} \theta_j^4 = R_1(\vec{\theta}).
	\end{equation}
	Using several times inequality (\ref{5.21}) for  $k=1,2,\ldots,n-1$ and combining with (\ref{5.22}), (\ref{5.17}), we obtain (\ref{P51}).  
	\hfill$\vspace{5mm}\scriptstyle\blacksquare$
	
\begin{Proof} {\it Lemma 4.1.} 	
Let
\begin{equation}
	F(\vec{\theta}) = \exp \left( -a\, \vec{\theta}^T \hat{Q} \vec{\theta} - b\sum\limits_{\{v_j,v_k\}\in EG} \Delta_{jk}^4 \right). 
\end{equation}
Note that for $\theta \in L$,
\begin{equation}
	F(\vec{\theta}) = \exp\left(
			 - a\sum\limits_{\{v_j,v_k\}\in EG}\Delta_{jk}^2
			- b\sum\limits_{\{v_j,v_k\}\in EG} \Delta_{jk}^4 
			\right).
\end{equation}
Since the integrand of (\ref{Eq_4_6}) is invariant under uniform
translation of all the
 $\theta_j$'s mod $\pi$ as well as  $V_0$ are mapped into itself, we can fix $\bar{\theta} = 0$ and multiply it by the ratio of its range $\pi$ to the length $n^{-1/2}$ of the vector $\frac{1}{n}[1,1,\ldots,1]^T$. 
	Thus we get that
		\begin{equation}\label{5.25}
			\begin{aligned}
		\int\limits_{V_0}	
			\exp\left(
			 - a\sum\limits_{\{v_j,v_k\}\in EG}\Delta_{jk}^2
			- b\sum\limits_{\{v_j,v_k\}\in EG} \Delta_{jk}^4 
			\right)
		 d\vec{\theta} = \\ 
 = \pi n^{1/2}\int\limits_{L\cap U_n(n^{-1/2+\varepsilon})}	
			F(\vec{\theta})
		 dL.
		 \end{aligned}
	\end{equation}
	In a similar way as in (\ref{LcapU_n->Un}) we find that
	\begin{equation}\label{5.26}
		\begin{aligned}
			\int\limits_{L \cap U_n(n^{-1/2+\varepsilon})}F(\vec{\theta})dL 
		 	 	= \int\limits_{P(\vec{\theta}) \in U_n(n^{-1/2+\varepsilon})} F(\vec{\theta})\, d \vec{\theta}
		 	\Big/  \int\limits_{-\infty}\limits^{+\infty} e^{-anx^2} dx = \\
		 	=\left(\frac{\pi}{an}\right)^{-1/2}\int\limits_{P(\vec{\theta}) \in U_n(n^{-1/2+\varepsilon})} F(\vec{\theta})\, d \vec{\theta}. 
		 	\end{aligned}
	\end{equation} 
	Using (\ref{1/2U_n}) and (\ref{U_n_in_P}), we get that
	\begin{equation}\label{5.27}
		\begin{aligned}
			\int\limits_{P(\vec{\theta}) \in U_n(n^{-1/2+\varepsilon})} F(\vec{\theta})\, d \vec{\theta} =
			 \int\limits_{U_n(n^{-1/2+\varepsilon})} F(\vec{\theta})\, d \vec{\theta} 
			+ O\left(\exp(-cn^{2\varepsilon})\right)\int\limits_{\mathbb{R}^n} e^{ - a\, \vec{\theta}^T \hat{Q} \vec{\theta} }.
		\end{aligned}
	\end{equation}
	Combining (\ref{P51}), (\ref{5.25}), (\ref{5.26}) and (\ref{5.27}), we obtain (\ref{Eq_4_6}).
\end{Proof}

\section{Final remarks}

Combining (4.2), (4.5), (4.6), (\ref{eq_Lemma_S1234}), (5.25), (5.26), (5.27), we find that:

\begin{Proposition} Let the assumtions of Theorem \ref{main} hold. Then 
\begin{equation}
	EO(G) = \left(1+ O\left(n^{-1+6\varepsilon}\right)\right) 2^{|EG|-1/2} \pi ^{-n+1/2} n \,Int,
\end{equation}
where 
\begin{equation}\label{888}
\begin{aligned}
Int = 
		\int\limits_{U_n(n^{-1/2+\varepsilon})}	
		\exp\left(
				 - \frac{1}{2}\, \vec{\theta}^T \hat{Q} \vec{\theta} 
			- \frac{1}{12}\sum\limits_{\{v_j,v_k\}\in EG} \Delta_{jk}^4 
			\right) 
		 d\vec{\theta},\\
	\hat{Q} = Q+J, \ \ \ \Delta_{jk} = \theta_j - \theta_k,		 	
\end{aligned}
\end{equation}
 where $Q=Q(G)$ denotes the Laplcian matrix and  $J$ is the matrix with every entry $1$.  
\end{Proposition}  
 Integral(\ref{888})
can be evaluated more precisely for specific classes of graphs in order to get asymptotic formulas 
for $EO(G)$ similar to (\ref{Eq_1_1}). For example, for even $n\rightarrow \infty$ 
\begin{equation}
		EO(K_{n,n}) = e^{-1}\frac{2^{n^2+n-\frac{1}{2}}}{\pi^{n-\frac{1}{2}} n^{n-1}} \left(1+ O(n^{-1/2+\varepsilon})\right)
	\end{equation} 
for any $\varepsilon>0$, where $K_{n,n}$ denotes the complete bipartite graph with $n$ vertices in each part.

\section*{Acknowledgements}
 This work was carried out under the supervision of S.P. Tarasov and supported in part by
RFBR grant no 11-01-00398a.

\section*{Warning}

{\bf It is preprint of the paper "Asymptotic behavior of the number of Eulerian orientations of graphs" 
accepted for publication in "Mathematical Notes" in 2012, (C) by Pleiades Publishing, Ltd., see also http://www.maik.ru.}

\noindent
{ {\bf M.I. Isaev}\\
Centre de Math\'ematiques Appliqu\'ees, Ecole Polytechnique,

91128 Palaiseau, France\\
Moscow Institute of Physics and Technology,

141700 Dolgoprudny, Russia\\
e-mail: \tt{isaev.m.i@gmail.com}}\\

\end{document}